\documentclass[12pt]{article}
\usepackage[T1]{fontenc}
\usepackage[english]{babel}
\usepackage{amsmath,latexsym,enumerate,amsthm,amssymb}

\newcommand{\al}{\alpha}

\newcommand{\de}{\delta}
\newcommand{\ep}{\epsilon}
\newcommand{\la}{\lambda}

\newcommand{\si}{\sigma}

\newcommand{\De}{\Delta}

\newcommand{\C}{\mathbb{C}}

\newcommand{\N}{\mathbb{N}}

\newcommand{\R}{\mathbb{R}}
\newcommand{\T}{\mathbb{T}}

\newcommand{\zerob}{\boldsymbol{0}}
\newcommand{\unb}{\boldsymbol{1}}

\newcommand{\kb}{\mathbf{k}}
\newcommand{\lb}{\boldsymbol{\ell}}
\newcommand{\mb}{\mathbf{m}}

\renewcommand{\sb}{\mathbf{s}}
\newcommand{\ub}{\mathbf{u}}
\newcommand{\vb}{\mathbf{v}}
\newcommand{\xb}{\mathbf{x}}
\newcommand{\yb}{\mathbf{y}}
\newcommand{\zb}{\mathbf{z}}

\newcommand{\Xb}{\mathbf{X}}

\newcommand{\ab}{\boldsymbol{\alpha}}
\newcommand{\bb}{\boldsymbol{\beta}}
\newcommand{\gb}{\boldsymbol{\gamma}}
\newcommand{\db}{\boldsymbol{\delta}}
\newcommand{\epb}{\boldsymbol{\epsilon}}
\newcommand{\mub}{\boldsymbol{\mu}}

\newcommand{\p}{\partial}

\newcommand{\ds}{\displaystyle}

\newcommand{\er}{{\cal E}}

\newtheorem{corollaire}{Corollary}

\newtheorem{definition}{Definition}

\newtheorem{rem}{Remark}

\newtheorem{theoreme}{Theorem}
\newtheorem*{theo A et B}{Theorem (\cite{crise})}

\newtheorem{notations/definitions}{Notations/Definitions}

\begin{document}
\title{Relations between values at $T$-tuples of negative integers \\
 of twisted multivariable zeta series \\
 associated to polynomials of several variables.}
\author{Marc de Crisenoy \footnote{Institut de 
Math\'ematiques de Bordeaux (A2X UMR 5465),
Universit\'e Bordeaux 1,
351 cours de la Lib\'eration
33405 Talence cedex, France.
Email: mdecrise@math.u-bordeaux1.fr
} 
 \hskip 0.5cm \& \hskip 0.5cm
Driss Essouabri\footnote{Universit\'e de Caen,
 UFR des Sciences (Campus II),
 Laboratoire de Math\'ematiques Nicolas Oresme (CNRS UMR 6139),
 Bd. Mal Juin, B.P. 5186, 14032 Caen, France.
Email: essoua@math.unicaen.fr}}

\maketitle \noindent {\small {{\bf {Abstract.}}
{\it We  give a new and very concise 
proof of  the existence 
of a  holomorphic continuation for a large class of twisted multivariable zeta functions. To do this,
we use a simple  method of ``decalage" that avoids using   an integral representation of the
zeta function. This allows us to derive explicit recurrence {\it relations} between the values at
$T-$tuples of negative integers. This also extends some earlier results of  several 
authors where the underlying polynomials were products of linear forms.} 


\section{Introduction}

\noindent Let $ Q,P_1,\dots,P_T \in \R[X_1,\ldots,X_N] $ and $
\mu_1,\ldots,\mu_N \in \C - \{1\},$ each   of modulus $ 1 $. To this data  we can
associate the following ``twisted" multivariable zeta series:
$$ Z(Q;P_1,\dots,P_T;\mu_1,\dots,\mu_N;s_1,\dots,s_T) = \sum_{m_1 \geq
1,\dots,m_N \geq 1} \frac{(\prod_{n=1}^N \mu_n^{m_n})
Q(m_1,\dots,m_N)}{\prod_{t=1}^T
P_t(m_1,\dots,m_N)^{s_t}}$$
where $ (s_1,\dots,s_T) \in \C^T $. \\
In this article we will always assume that:
$$ \forall t \in \{1,\dots,T \} \ \ \forall \xb \in [1, + \infty[^N \
P_t(\xb) > 0 {\mbox { and }}   
\prod_{t=1}^T P_t(\xb) \xrightarrow[\substack{ |\xb| \rightarrow
+\infty \\ \xb \in J^N}]{} + \infty  \qquad (\#)$$
It is not difficult to see that {\bf condition (\#)} implies   that  $
Z(Q;P_1,\dots,P_T;\mu_1,\dots,\mu_N;s_1,\dots,s_T)
$ is an absolutely convergent series when  $ \Re(s_1),\dots,\Re(s_T) $ are
sufficiently large. \\
Cassou-Noguès (\cite{cassou 82}) and  Chen-Eie (\cite{chen; eie})
prooved in the case $ T=1, $ and   $ P = P_1 $ a polynomial with
positive coefficients, that the above series can be {\it holomorphically}
continued to the whole complex plane and obtained very nice formulas for
their values at negative integers.
In \cite{crise} de Crisenoy extended these results by allowing $T > 1$ and introducing
the (HDF) hypothesis (see definition \ref{hdf} in \S 2) that is much
weaker than positivity of coefficients. 
The Main result of \cite{crise} is the following:

\begin{theo A et B}
Let $ Q, P_1,...,P_T \in \R[X_1,...,X_N] $ and $ \mub \in (\T \setminus \{1\})^N $. \\
Assume that: \\
For each $t=1,\dots,T$ $P_t$  satisfies the (HDF) hypothesis and 
that $ \ds \prod_{t=1}^T P_t(\xb) \xrightarrow[ \substack{ |\xb| \rightarrow + \infty \\ \xb \in J^N}]{} + \infty $. \\
Then: 
\begin{itemize}
\item $ Z(Q;P_1,\dots,P_T;\mub;\cdot) $ can be holomorphically extended to $
\C^T $;
\item 
 For all  $ k_1,...,k_T \in \N $ if we set $ \ds Q
\prod_{t=1}^T P_t^{k_t}  = \sum_{\ab \in S} a_{\ab} \Xb^{\ab} $, we
have:
$$ Z(Q;P_1,...,P_T;\mub;- k_1,...,- k_T) = \sum_{\ab \in S} a_{\ab}
\prod_{n=1}^N \zeta_{\mu_n}(-\al_n).$$
where for all $\mu \in \T$, $\zeta_{\mu}(s)=\sum_{m=1}^{+\infty}
\frac{\mu^m}{m^s}$.
\end{itemize}
\end{theo A et B}
To obtain this theorem he used an integral 
representation for the zeta series. The proof of the holomorphic
continuation of the resulting integral is long and complicated.\\ 
By restricting to hypoelliptics polynomials (see definition
\ref{hypoelliptic} in \S 2) we can avoid the integrals:\\
the first result of this article gives  a new proof of the 
{\it holomorphic} continuation of these series  
under the assumption that $ P_1,\dots,P_T $ are hypoelliptic. Our
proof uses the ''decalage'' method of Essouabri (\cite{essouabri}). 
Of particular interest for this paper is that this method does {\it not} use an
integral representation for the zeta series. And the resulting proof
is very concise and much simple.\\
The second result is relations between the values at $T-$tuples of negative
integers of these series, relations that are true under the $ HDF $
hypothesis, and that give a mean to calculate by induction the values. \\
These relations are very simple in the case of linear forms, particularly
interesting because of its link with the zeta functions of number fields
(\cite{cassou 79}). Similar results have been obtain by severals
authors in particular cases of linear forms (see  \cite{akiyama; egami; tanigawa}, 
\cite{akiyama; tanigawa} and \cite{arakawa; kaneko}). 
Our method allows one also to obtain new relations even in the cases
of linear forms. 

\section{Notations and preliminaries:}
First, some notations:
\begin{enumerate}
\item Set $ \N = \{0, 1, 2, \dots\}, \N^* = \N - \{0\}, \  J = [1, +
\infty[ ,\ $ and $ \T =\{ z \in \C\ | \ |z| = 1 \} $. 
\item 
The real part of $ s \in \C $ will be denoted $ \Re(s) = \si $ and its
imaginary part $ \Im(s) = \tau $. 
\item 
Set $ \zerob = (0,\ldots,0) \in \R^N $ and $ \unb = (1,\ldots,1) \in \R^N
$. 
\item
For $ \xb = (x_1,\ldots,x_N) \in \R^N $ we set $ |\xb| = |x_1|+\ldots+|x_N|
$. 
\item
For $ \zb = (z_1,\ldots,z_N) \in \C^N $ and $ \ab = (\al_1,\ldots,\al_N)
\in \R_+^N $ we set $ \zb^{\ab} = z_1^{\al_1} \ldots z_N^{\al_N} $. 
\item 
The notation $ f(\la,\yb,\xb) \ll_{\yb} g(\xb) \ \ (\mbox
{uniformly in } 
\xb \in X \mbox  { and } \la \in \Lambda) $ means that there exists $ A =
A(\yb) > 0 $, that depends neither on $ \xb $ nor
$ \la $, but   could a priori depend  on other parameters, and in
particular on $ \yb $, such that: $ \forall \xb \in X \ \forall \la \in
\Lambda \ |f(\la,\yb,\xb)| \leq A g(\xb) $. \\
When there is no ambiguity, we will omit the word uniformly and the index $
\yb. $  
\item
The notation $ f \asymp g $ means that we  have both $ f \ll g $ and $ g
\ll f $. 
\item 
For ${\bf a} \in \R_+^N$ and  $ P \in \R[X_1,\dots,X_N] $, we define $ \De_{\bf a}  P \in \R[X_1,\dots,X_N] $
by\\ $ \De_{\bf a} P(\Xb) = P(\Xb + {\bf a}) - P(\Xb) $. If no ambiguity, we
will note simply $ \De_{\bf a} P$ by $\De P$.
\item
Let ${\bf a} \in \N^N$, $N \in \N^*$, $I \subset \{1,\dots,N\}$,
$q=\#I$, $I^c= \{1,\dots,N\}\setminus I$ and \\
${\bf b}= (b_i)_{i\in I^c} \in \prod_{i\in I^c} \{0,\dots,a_i\}$.
Set $I=\{i_1,\dots,i_q\}$ and $I^c=\{j_1,\dots,j_{N-q}\}$. Then:
\begin{enumerate}
\item For all $ {\bf d}= (d_1,\dots,d_q) \in \N^q $ we define  
${\bf m}^{{\bf a},I,{\bf b}}({\bf d})=(m_1,\dots,m_N)\in \N^N$ by \\
$\forall k=1,\dots,q$, $m_{i_k}= a_{i_k}+d_k$ and $\forall
k=1,\dots,N-q$, $m_{j_k}=b_{j_k}$;
\item For all $f :\N^N \rightarrow \C$ we define 
$f^{{\bf a},I,{\bf b}}:\N^q \rightarrow \C$ by 
$f^{{\bf a},I,{\bf b}}({\bf d})=
f\left({\bf m}^{{\bf a},I,{\bf b}}({\bf d})\right)$.
\end{enumerate}
\end{enumerate}

\noindent {\bf Convention}: In this work we will say that a series defined
by a sum over $ N \geq 1 $ variables is 
convergent when it is absolutely  convergent. \\ \\

Let us recall some definitions: \\

\begin{definition} \label{hypoelliptic} 
 $ P \in \R[X_1,\ldots,X_N] $ is said to be hypoelliptic if: \\
$ \forall \xb \in J^N \ P(\xb) > 0 $ \ and \ $ \ds \forall \ab \in \N^N
\setminus \{ \zerob \} \ \ \frac{\p^{\ab} P}{P}(\xb)
\xrightarrow[\substack{|\xb| \rightarrow
+\infty \\ \xb \in J^N}]{} 0 $. 
\end{definition}

\begin{definition} \label{hdf}
Let $ P \in \R[X_1,\ldots,X_N] $. \\ 
$ P $ is said to satisfy the weak decreasing hypothesis (denoted  $HDF$ in
the rest of the article) if: 

\begin{itemize}
\item $ \forall \xb \in J^N \ P(\xb) > 0 $, 
\item $ \exists \ep_0 > 0 $ such that for $ \ab \in \N^N $ and $ n \in \{
1, \dots ,N \} $: $ \ds \al_n \geq 1 \Rightarrow \frac{{\p}^{\ab} P
}{P}(\xb) \ll x_n^{-\ep_0} \ \ ( \xb \in J^N) $.
\end{itemize}
\end{definition}
\begin{rem}
It follows from H{\"{o}}rmander (\cite{hormander}, p.62) (see also Lichtin
(\cite{lichtin}, P.342) that if \\
$ P \in \R[X_1,\ldots,X_N] $ is hypoelliptic then there exists 
$ \ds \ep > 0 $ such that $ \ds \forall \ab \in \N^N \setminus \{
\zerob \}$\\ $\ds \frac{ \p^{\ab} P }{P} (\xb) \ll {\xb}^{-\ep \unb} \ (\xb \in J^N) $. 
Therefore $P$ satisfy also the HDF hypothesis.
\end{rem}

\section{Main results:}

We first give a new and concise proof of the existence of a
holomorphic continuation for our series when each polynomial is
hypoelliptic. The procedure uses the ''decalage'' method, 
and does not require an integral representation.\\
A consequence of this method is found in corollaries 1 and 2. 
This gives simple recurrent {\it relations} between the values at
$T-$tuples of negative integers when the polynomials $P_i$ are
products of linear polynomials or a class of quadratic polynomials. 
This procedure is different from that in \cite {crise} where   
the calculation  of the special values was not an immediate 
consequence of the holomorphic continuation of the zeta series. \\

\begin{theoreme}
Let $ \mub \in (\T \setminus \{1 \} )^N $ and $ Q,P_1,\dots,P_T \in
\R[X_1,\dots,X_N] $. \\
We assume that $ P_1,\dots,P_T $ are hypoelliptic and that at least 
one of them is not constant. \\
Then $ Z(Q;P_1,\dots,P_T;\mub;\cdot) $ can be holomorphically extended 
to $ \C^T $. \\
Let $ {\bf a} \in \N^N $. We set $ \ds  E({\bf a}) = \{ \mb\in \N^{*N} \ |
\ \mb  \ngeq  {\bf a} + \unb \} = \left\{ \mb \in \N^{*N} \ | \ \exists n
\in \{ 1,\dots,N \} \ m_n \leq a_n \right \} $. \\
Then for all $ \kb \in \N^T $ we have the following relation:
\begin{align*}
(1- \mub^{{\bf a}}) Z(Q;P_1,\dots,P_T;\mub;-\kb) = & \mub^{{\bf a}}  \sum_{
\zerob \leq \ub < \kb }  \binom{\kb}{\ub}  Z \left( Q(\Xb + {\bf a})
\prod_{t=1}^T  (\De_{\bf a} P_t)^{k_t-u_t};P_1,\dots,P_T; \mub; -\ub \right) \\
                                   & + \mub^{{\bf a}} 
Z(\De Q ; P_1,\dots,P_T; \mub ;-\kb) + Z_{N-1}^{{\bf a}}(-\kb).
\end{align*} 
Where $\ds Z_{N-1}^{\bf a}(\sb) =  \sum_{ \mb \in E({\bf a}) } \mub^{\mb}
Q(\mb) \prod_{t=1}^T P_t  (\mb)^{- s_t}$.
By using notation (9) above it's easy to see that:
$$Z_{N-1}^{\bf a}(\sb) = \sum_{q=1}^{N-1} \sum_{I\subset \{1,\dots,N\} \atop \# I =q }
\sum_{{\bf b}=(b_i) \in \prod_{i\in I^c} \{0,\dots,a_i\}}
\left(\prod_{i \in I^c} \mu_i^{b_i}\right)~Z\left(Q^{{\bf a},I,{\bf b}};P_1^{{\bf a},I,{\bf b}},\dots,
P_T^{{\bf a},I,{\bf b}};(\mu_i)_{i\in I};\sb \right).$$
In particular, it is clear that $Z_{N-1}^{\bf a} $ is a finite
linear combination of zeta series $ Z $ associated to hypoelliptic
polynomials of at most $ N - 1 $ variables. 
\end{theoreme}

\begin{rem} 
When $ \mub^{{\bf a}} \neq 1 $ (and we can of course choose $ {\bf a} $
such that this is satisfied, this formula allows us to compute the values
at $T-$tuples of negative integers of the series $ Z $ by recurrence
because in each term on the right a integral value is strictly less than
the corresponding value on the left: \\
$ |\ub| < |\kb| $, $ \deg(\De Q) < \deg Q $, $ N-1 < N $. 
\end{rem}

Now, using theorem A of \cite{crise}(that gives the existence of the
holomorphic continuation under the HDF hypothesis) and theorem B of
\cite{crise} (that gives closed formula for the values, still under HDF)
and the preeceding theorem, we show that the relations remain true under
the HDF hypothesis: \\

\begin{theoreme}
Let $ \mub \in (\T \setminus \{1 \} )^N $ and $ Q,P_1,\dots,P_T \in
\R[X_1,\dots,X_N] $. \\
We assume that $ P_1,\dots,P_T $ satisfies the HDF hypothesis and that
$ \ds \prod_{t=1}^T P_t(\xb) \xrightarrow[\substack{ |\xb| \rightarrow
+\infty \\ \xb \in J^N}]{} + \infty $.\\
Then, the relations of theorem 1 are still true.
\end{theoreme}

Now we deal with the particular case of linear forms. In this case the
relations become particularly simple! \\ 

\begin{corollaire}
Let $ \mub \in (\T \setminus \{1 \} )^N $ and $ L_1,\dots,L_T $ linear
forms with positive coefficients. \\
We assume that for each $ n $ there exists  $ t $ such that $ L_t $ really
depends on $ X_n $. \\
We set $ Z(\mub;\cdot) =  Z(1;L_1,\dots,L_T;\mub;\cdot) $. \\
Let $ {\bf a} \in \N^N $. We set $ \ds  E({\bf a}) = \{ \mb\in \N^{*N} \ |
\ \mb  \ngeq  {\bf a} + \unb \} = \left\{ \mb \in \N^{*N} \ | \ \exists n
\in \{ 1,\dots,N \} \ m_n \leq a_n \right \} $. \\
For all $ t \in \{1,\dots,T \} $ we set: $ \de_t = L_t(\Xb +{\bf a}) -
L_t(\Xb) $. $ \de_t \in \R $. \\
Then, for all $ \kb \in \N^T $ we have the following relation:
$$ (1- \mub^{{\bf a}}) Z( \mub;-\kb) =  \mub^{{\bf a}}  \sum_{ \zerob \leq \ub
< \kb } \db^{\kb-\ub} \binom{\kb}{\ub}  Z \left(\mub; -\ub \right) +
Z_{N-1}^{{\bf a}}( \mub,-\kb) $$
Where $\ds Z_{N-1}^{\bf a}(\sb) =  \sum_{ \mb \in E({\bf a}) } \mub^{\mb} 
\prod_{t=1}^T L_t  (\mb)^{- s_t}$.
By using notation (9) above it's easy to see that:
$$
Z_{N-1}^{\bf a}(\sb) = \sum_{q=1}^{N-1} \sum_{I\subset \{1,\dots,N\} \atop \# I =q }
\sum_{{\bf b}=(b_i) \in \prod_{i\in I^c} \{0,\dots,a_i\}}
\left(\prod_{i \in I^c} \mu_i^{b_i}\right)~Z\left(1;L_1^{{\bf a},I,{\bf b}},\dots,
L_T^{{\bf a},I,{\bf b}};(\mu_i)_{i\in I};\sb \right).$$
In particular, it is clear that $Z_{N-1}^{\bf a} $ is a finite
linear combination of zeta series $ Z $ associated to linear
forms of at most $ N - 1 $ variables. 
\end{corollaire}

\begin{rem}
Of course, here as well this formulaes allows a calculus by induction. \\
\end{rem}

Since the works of Cassou-Noguès and Shintani, we know that the case of
linear forms is particulaly interesting for algebraic number theory because
of the link with the zetas functions of number fields. See also
\cite{zagier 94} for more motivations.

\begin{rem}
Let us assume that $ P \in \R[X_1,\dots,X_N] $ is a product of linear
forms: $ \ds P = \prod_{t=1}^T  L_t $ where $ L_1,\dots,L_T $ have real
positive coefficients, and that we want to evaluate the numbers $
Z(1;P;\mub;-k) $ where $ k \in \N $. \\
We could use theorem I, but it seems more interesting to note that \\ $
\forall s \in \C \ \  Z(1;P;\mub;s) = Z(1;L_1, \dots,L_T;\mub; s, \ldots, s)
$ and then to consider the numbers \\
$ Z(1;P;\mub;-k) = Z(1;L_1, \dots,L_T;\mub;-k, \ldots,-k) $ inside the
family $ Z(1;L_1, \dots,L_T;\mub;-k_1, \ldots,-k_T) $ because of the simple
relations between these numbers given by the preceeding proposition.
\end{rem}

The particular case of linear forms is not the only case when 
relations  of theorems 1 and 2 become particularly simple. The
great flexibility in the  assumptions of this theorems, allowed one to
obtain also very simple relations  in somes others cases as in the
following:

\begin{corollaire}
Let $ \mub \in (\T \setminus \{1 \} )^N $ and let $ P_1,\dots,P_T \in
\R[X_1,\dots,X_N]$ be polynomials of degree at most $2$.
Suppose that for all $t=1,\dots,T$:
$$P_t({\bf X})=P_t(X_1,\dots,X_N)=\sum_{k=1}^{r_t} \left(\langle \alpha^{t,k} ,
  {\bf X}\rangle \right)^2 + \sum_{n=1}^N c_{t,n} X_n+d_t$$
where $\alpha^{t,k} \in \R^N$, $c_{t,i}\in \R_+^*$ and $d_t \in \R_+$
 $\forall t,k,i$.\\
Assume that there exist ${\bf a} \in \N^N \setminus \{\zerob \}$ 
such that $\langle \alpha^{t,k} ,{\bf a} \rangle =0$ for all $t,k$.\\
Then:
\begin{enumerate}
\item for all $ t \in \{1,\dots,T \} $,  $ \de_t := P_t(\Xb +{\bf a}) -
P_t(\Xb) \in \R$; 
\item the relations of Corollary 1 are still true.
\end{enumerate}
\end{corollaire}

\section{Proof of theorem 1} 
Let $t\in \{1,\dots,T\}$. It's clear that there exists $\ab \in
\N^N$ such that $\partial^{\ab} P_t$ is a no vanishing constant
polynomial. So the hypoellipticity of $P_t$ implies that 
$P_t({\bf x})\rightarrow +\infty$ when $|{\bf x}|
\rightarrow +\infty$ $({\bf x} \in [1,+\infty[^N )$.
By using Tarski-Seidenberg (see for example \cite{essouabri1}, Lemme 1), 
we know that there exists $\delta >0$ such 
that $P_t({\bf x})\gg (x_1\dots x_N)^\delta$
uniformly in ${\bf x} \in [1,+\infty[^N$.
Therefore, its clear that there exists $ \si_0 $ such that if $ \si_1,\dots,\si_T > \si_0 $ then $
Z(Q,P_1,\dots,P_T,\mub,\sb) $ converges. \\
By Remark 2.1 above, there exist also  a fixed $ \ds \ep > 0 $ such that $ \forall t \in \{ 1,\dots,T \} $ we
have $$ \ds \forall \ab \in \N^N \setminus \{ \zerob \} \quad 
\frac{ \p^{\ab} P}{P} (\xb) \ll {\xb}^{-\ep \unb} \ (\xb \in J^N).$$

\noindent {\bf Step 1}: we establish a formula $(*)$. \\

\noindent {\bf Proof of step 1}: \\
$ P_1,\dots,P_T $ are fixed in the whole proof so we will denote $
Z(Q,\mub,\cdot) $ instead of $ Z(Q,P_1,\dots,P_T,\mub,\cdot) $. \\
With this notation for all $ \sb $ such that $ \si_1,\dots,\si_T > \si_0 $
we have:
\begin{align*}
 Z(Q,\mub,\sb) & = \sum_{\mb \in \N^{*N}} \mub^{\mb} Q(\mb) \prod_{t=1}^T
P_t(\mb)^{- s_t} \\
          &  =  \sum_{\mb  \geq {\bf a} + \unb} \mub^{\mb} Q(\mb)
\prod_{t=1}^T P_t  (\mb)^{- s_t}  + Z_{N-1}^{{\bf a}}(\sb) \\
          &  =  \mub^{{\bf a}}  \sum_{\mb \in \N^{*N}}  \mub^{\mb} Q(\mb +
{\bf a}) \prod_{t=1}^T P_t  (\mb + {\bf a})^{- s_t}  + Z_{N-1}^{{\bf a}}(\sb)
\end{align*}

\noindent $ U \in \N $ is set for the whole step. \\
We have $ g_U \colon \C \times \C \setminus ]-\infty,-1] \to \C $
holomorphic and satisfying: \\
$ \ds \forall s \in \C  \ \forall z \in \C \setminus ]-\infty, -1] \ \
(1+z)^s = \sum_{u=0}^U \binom{s}{u} z^u+z^{U+1} g_U(s,z) $. \\
$ \forall k \in \N $ vérifiant $ k \leq U $ et $ \forall z \in \C \setminus
]-\infty,-1] $, $ g_U (k,z) = 0 $. \\

For $ t \in \{1,\dots,T \} $ we define $ \De_t = \De_{\bf a} P_t $. \\
\noindent For $ t \in \{ 1,\dots,T \} $ and $ \mb \in \N^{*N} $, we define
$ H_{t,\mb,U} \colon \C \to \C $ by: 
$$  H_{t,\mb,U}(s_t) = \sum_{u_t =
0}^U \binom{-s_t}{u_t} \De_t(\mb)^{u_t} P_t(\mb)^{-u_t}.$$

\noindent $ \forall t \in \{1,\dots,T \} $ we have:
\begin{align*}
 P_t(\mb + {\bf a})^{- s_t} & = \left[P_t(\mb) + \De_t(\mb) \right]^{- s_t} \\
                          & =  P_t(\mb)^{- s_t} \left[ 1 + \De_t(\mb)
P_t(\mb)^{-1} \right]^{- s_t} \\
                          & =  P_t(\mb)^{- s_t}  \left[ H_{t,\mb,U}(s_t) +
\De_t(\mb)^{U+1} P_t(\mb)^{-(U+1)} g_U \left(- s_t, \De_t(\mb)
P_t(\mb)^{-1} \right) \right] 
\end{align*}

\noindent For $ x_1,\dots,x_T,y_1,\dots,y_T \in \R $, we have $ \ds
\prod_{t=1}^T (x_t+y_t) = \sum_{\epb \in \{0,1 \}^T}  \prod_{t=1}^T x_t^{1-
\ep_t} y_t^{\ep_t} $, so:
$$ \prod_{t=1}^T  P_t(\mb + {\bf a})^{- s_t} =  \sum_{\epb \in \{0,1 \}^T}
\prod_{t=1}^T  H_{t,\mb,U}(s_t)^{1- \ep_t} \De_t(\mb)^{\ep_t (U+1)}
P_t(\mb)^{-s_t - \ep_t (U+1)} g_U \left(- s_t, \De_t(\mb) P_t(\mb)^{-1}
\right) ^{ \ep_t} $$

For $ \mb \in \N^{*N} $ and $ \epb \in  \{0,1 \}^T $ we define $
f_{\mb,U,\epb} \colon \C^T \to \C $ thanks to the following formula:
$$ f_{\mb,U,\epb}(\sb) =  \prod_{t=1}^T  H_{t,\mb,U}(s_t)^{1- \ep_t}
\De_t(\mb)^{\ep_t (U+1)} P_t(\mb)^{-s_t - \ep_t (U+1)} g_U \left(- s_t,
\De_t(\mb) P_t(\mb)^{-1} \right) ^{ \ep_t} $$
So for all $ \mb \in \N^{*N} $ and $ \sb \in \C^T $ we have:
$ \ds \prod_{t=1}^T  P_t(\mb + {\bf a})^{- s_t} =  \sum_{\epb \in \{0,1
\}^T} f_{\mb,U,\epb} (\sb) $. \\
We define $ Z_U (Q,\mub,\cdot) $ by: $ \ds Z_U (Q,\mub,\sb) :=  \sum_{\epb
\in \{0,1 \}^T \setminus \{ \zerob \} } \  \sum_{\mb \in \N^{*N}}
\mub^{\mb} Q(\mb + {\bf a})  f_{\mb,U,\epb}(\sb) $. \\
We will see in step 2 that for $ U $ large enough $ Z_U (Q,\mub,\cdot) $
exists and is holomorphic on 
$$ \left\{ \sb \in \C^T \ | \ \si_1,\dots,\si_T
> \si_0 \right\}.$$ 
It's clear that 
$ \ds Z_U(Q,\mub,\sb) = \sum_{\mb \in \N^{*N}} \mub^{\mb} Q(\mb + {\bf a})
\sum_{\epb \in \{0,1 \}^T \setminus \{ \zerob \} } f_{\mb,U,\epb}(\sb).$\\
So $ \ds Z(Q,\mub,\sb) = \mub^{{\bf a}} \sum_{\mb \in \N^{*N}} \mub^{\mb}
Q(\mb + {\bf a}) f_{\mb,U,\zerob}(\sb) + \mub^{{\bf a}} Z_U (Q,\mub,\sb) +
Z_{N-1}^{{\bf a}}(\sb) $. 

\noindent By definition $ \ds f_{\mb,U,\zerob}(\sb) = \prod_{t=1}^T
H_{t,\mb,U}(s_t) P_t(\mb)^{-s_t} $ so:
\begin{align*}
f_{\mb,U,\zerob}(\sb) & = \prod_{t=1}^T \sum_{u_t = 0}^U \binom{-s_t}{u_t}
\De_t(\mb)^{u_t} P_t(\mb)^{-(s_t + u_t)} \\
                      & = \sum_{0 \leq u_1,\dots,u_T \leq U} \prod_{t=1}^T
\binom{-s_t}{u_t} \De_t(\mb)^{u_t} P_t(\mb)^{-(s_t + u_t)} \\
                      & = \sum_{ \ub \in \{0,\dots,U \}^T}
\binom{-\sb}{\ub} \prod_{t=1}^T  \De_t(\mb)^{u_t} P_t(\mb)^{-(s_t + u_t)} 
\end{align*}
then, for all $ \sb $ such that $ \si_1,\dots,\si_T > \si_0 $, we have:
\begin{align*}
 \sum_{\mb \in \N^{*N}} \mub^{\mb} Q(\mb + {\bf a}) f_{\mb,U,\zerob}(\sb) &
= \sum_{\mb \in \N^{*N}} \mub^{\mb} Q(\mb + {\bf a})  \sum_{ \ub \in
\{0,\dots,U \}^T} \binom{-\sb}{\ub}  \prod_{t=1}^T  \De_t(\mb)^{u_t}
P_t(\mb)^{-(s_t + u_t)} \\
                                                                      & =
\sum_{ \ub \in \{0,\dots,U \}^T} \binom{-\sb}{\ub} \sum_{\mb \in \N^{*N}}
\mub^{\mb} Q(\mb + {\bf a}) \prod_{t=1}^T  \De_t(\mb)^{u_t} P_t(\mb)^{-(s_t
+ u_t)} \\ 
                                                                      & =
\sum_{ \ub \in \{0,\dots,U \}^T} \binom{-\sb}{\ub} Z \left( Q(\Xb + {\bf
a}) \prod_{t=1}^T  \De_t^{u_t} ,\mub, \sb + \ub \right) 
\end{align*}
The combination of the precedings results gives us:
$$ Z(Q,\mub,\sb) = \mub^{{\bf a}} \sum_{ \ub \in \{0,\dots,U \}^T}
\binom{-\sb}{\ub} Z \left( Q(\Xb + {\bf a}) \prod_{t=1}^T  \De_t^{u_t}
,\mub, \sb + \ub \right)  + \mub^{{\bf a}} Z_U (Q,\mub,\sb) + Z_{N-1}^{{\bf
a}}(\sb) $$
The following $(*)$ formula is now clear:
\begin{align*}
(*) \ \ \ \ \ \ \ \ \ \ \ (1- \mub^{{\bf a}}) Z(Q,\mub,\sb) = &  \mub^{{\bf
a}}  \sum_{ \ub \in \{0,\dots,U \}^T \setminus \{ \zerob \} }
\binom{-\sb}{\ub}  Z \left( Q(\Xb + {\bf a}) \prod_{t=1}^T  \De_t^{u_t}
,\mub, \sb + \ub \right) \\
                                  & + \mub^{{\bf a}} Z(\De Q,\mub,\sb) +
Z_{N-1}^{{\bf a}}(\sb) + \mub^{{\bf a}} Z_U(Q,\mub,\sb) 
\end{align*}

\noindent {\bf Step 2}: For all $ a \in \R_+ $ there exists $ U_0 \in \N $
such that for all $U\geq U_0$, $ Z_U (Q, \mub,\cdot) $ exists 
and is holomorphic on $ \left\{
\sb \in \C^T \ | \ \forall t \in \{1,\dots,T\} \ \si_t > -a \right\} $. \\

\noindent {\bf Proof of step 2}: \\
Let $ a \in \R $. Let $ K $ be a compact of $ \C^T $ included in $ \ds
\left\{ \sb \in \C^T \ | \ \forall t \in \{1,\dots,T \} \ \si_t > -a
\right\} $. \\ 
$ \star $ Let $ t \in \{1, \dots, T \} $. \\
We know that $ P_t(\xb) \gg 1 \ \ ( \xb \in  J^N ) $ so it's easy to
seen that : 
$ P_t(\xb)^{-s_t} \ll P_t(\xb)^a \ \ ( \xb \in J^N  \ \sb \in K) $. 

\noindent Let denote $ p = \max \left\{ \deg_{X_n} P_t \ | \ 1 \leq n \leq
N \ 1 \leq t \leq T \right\} $. \\
From now we assume $ a > 0 $. So $ \ds P_t(\xb)^a \ll \xb^{pa \unb} \ \
( \xb \in J^N ) $. \\ 
As a conclusion we have: $ \ds P_t(\xb)^{-s_t} \ll \xb^{pa \unb} \ \ ( \xb
\in J^N \ \ \sb \in K) $. \\

\noindent Let $ U \in \N^* $ and $ \epb \in \{0,1 \}^T $. By definition $
\ds H_{t,\mb,U}(\sb) = \sum_{u_t = 0}^U \binom{-s_t}{u_t} \left(
\frac{\De_t(\mb)}{P_t(\mb)} \right)^{u_t} $, so hypoellipticity of $
P_t $ implies that $ H_{t,\mb,U}(\sb) \ll 1 \ ( \mb \in \N^{*N}, \sb \in K)
$. \\
It implies also that there exists a compact of $ ]-1, +\infty [$ 
containing all the  $ \ds
\frac{\De_t(\mb)}{ P_t(\mb)} $ where $ \mb $ is in  $ \N^{*N} $ so: 
$ \ds g_U \left(- s_t, \De_t(\mb) P_t(\mb)^{-1} \right) ^{ \ep_t} \ll 1 \ (
\mb \in \N^{*N} \ \ \sb \in K) $. \\
Thanks to the Taylor formula and to the choice of $ \ep $ we have $ \ds
\frac{\De_t(\mb)}{P_t(\mb)} \ll \mb^{- \ep \unb} \ (\mb \in \N^{*N}) $. \\
From what preceeds we deduce that we have, uniformely in $ \mb \in \N^{*N}
$ and $ \sb \in K $: \\
$ \ds H_{t,\mb,U}(\sb)^{1- \ep_t} \left( \frac{\De_t(\mb)}{P_t(\mb)}
\right)^{\ep_t (U+1)} P_t(\mb)^{-s_t} g_U \left(- s_t, \De_t(\mb)
P_t(\mb)^{-1} \right) ^{ \ep_t} \ll \mb^{pa \unb} \mb^{- \ep \ep_t (U+1)
\unb} $. \\

\noindent $ \star $ By definition 
$$ f_{\mb,U,\epb}(\sb) = \prod_{t=1}^T  H_{t,\mb,U}(\sb)^{1- \ep_t} \left(
\frac{\De_t(\mb)}{P_t(\mb)} \right)^{\ep_t (U+1)} P_t(\mb)^{-s_t} g_U
\left(- s_t, \De_t(\mb) P_t(\mb)^{-1} \right) ^{ \ep_t} $$

Since $ \epb \neq \zerob $, we have: $ f_{\mb,U,\epb}(\sb) \ll \mb^{Tpa
\unb} \mb^{- \ep (U+1) \unb} \ \ (\mb \in \N^{*N}) $. \\
We denote $ q = \max \left\{ \deg_{X_n} Q \ | \ 1 \leq n \leq N \right\} $
(obviously we can assume that $ Q \neq 0 $). \\
We see that $ Q(\mb + {\bf a}) f_{\mb,U,\epb}(\sb) \ll \mb^{(q + Tpa - \ep
(U+1)) \unb} \ (\mb \in \N^{*N}) $. \\
So it is enough to have $ q + Tpa - \ep (U+1) \leq -2 $, for $ U $ to fit. \\
Thus it is enough to choose $ U \in \N $ such that $ \ds U \geq 
U_0:=\left[\frac{q + Tpa+2}{\ep}\right]+1$. \\

\noindent {\bf Convention}: \\
let us take a convention, that we will use until the end of the proof.
Let $ a \in \R $. We will say that a function $ Y $ is an entire
combination until $ a $ of the functions $ Y_1,\dots,Y_k $ if there exists: \\
$ \star $  entire functions $ \ds \la_1,\dots,\la_k \colon \C^T \to \C $, \\
$ \star $ one function $ \la \colon \left\{ \sb \in \C^T \ | \ \forall t
\in \{1,\dots,T \} \ \si_t > a \right\} \to \C $ holomorphic, \\
such that $ \ds Y = \la + \sum_{i=1}^k \la_i Y_i $. \\

\noindent {\bf A definition and a remark}: \\
for $ \ub \in \N^T $ and $ Q \in \R[X_1,\dots,X_N] $, we denote $ \er_{\ub}
(Q) $ the subspace of $ \R[X_1,\dots,X_N] $ generated by the polynomials of
the following form: $ \ds \p^{\bb} Q \prod_{t=1}^T \prod_{k \in F_t}
\p^{f_t(k)} P_t $ where: 

$ \star \ \bb \in \N^N $, 

$ \star $ $ F_1,\dots,F_T $ are finite subset of $ \N $  satisfying $
|F_t|=u_t $, 

$ \star \ \forall t \in \{1,\dots,T \} \ \ f_t $ is a function from $ F_t $
into $ \N^N \setminus \{ \zerob \} $. \\

\noindent We remark that $ \er_{\ub} (Q) $ is stable under partial
derivations. \\

We are now beginning the proof of the existence of an holomorphic
continuation. \\
The proof is by recurrence on $ N $; it will be clear that the proof that
rank $ N-1 $ implies rank $ N-1 $ gives the result at rank $ N = 1 $. \\
Let $ N \geq 1 $. We assume that the result is true at rank $ N-1 $. \\

\noindent {\bf Step 3}: \\
Let $ Q \in \R[X_1,\dots,X_N] $ and $ a \in \R_+ $. \\
Then $ Z(Q,\mub,\sb) $ is an entire combination until $ -a $ of functions
of the type $ Z(R,\mub,\sb+\ub) $ where $ \ub \in \N^N \setminus \{\zerob
\} $ and $ R \in \er_{\ub} (Q) $. \\

\noindent {\bf Proof of step 3}: \\
we are going to show by recurrence on $ d \in \N $ that if $ \deg Q < d $
then the result is true. \\
For $ d = 0 $ it is clear. \\
Let us assume the result for $ d \geq 1 $. \\
Let $ Q \in \R[X_1,\dots,X_N] $ such that $ \deg Q < d+1 $. \\
Thanks to step 2 we set $ U $ such that $ Z_U(Q,\mub,\cdot) $ is
holomorphic on $ \left\{ \sb \in \C^T \ | \ \forall t \in \{1,\dots,T \}
\ \si_t > -a \right\} $. \\
We are now going to use the formula $(*)$ obtained in step 1. We look at each
of the 4 term on the right.\\
$ \bullet $ Let $ \ub \in \{0,\dots,U \}^T \setminus \{ \zerob \} $. \\
Thanks to the Taylor formula, it is easy to see that $ \ds Q(\Xb + {\bf a})
\prod_{t=1}^T \De_t^{u_t} \in \er_{\ub}(Q) $. \\
$ \bullet $ $ \deg \De Q < d $ so the recurrence hypothesis
on $d$ implies that, $Z(\De Q, \mub,\sb) $ is an entire combination 
until $ -a $ of functions of the type $ Z(R,\mub,\sb + \ub) $ where $ \ub \in \N^N \setminus \{ \zerob
\} $ and $ R \in \er_{\ub} (\Delta Q) $. \\
Furthermore, clearly, $ \er_{\ub} (\Delta Q) \subset \er_{\ub}(Q) $. \\
$ \bullet $ Thanks to the recurrence hypothesis on $N$, $ Z^{\bf a}_{N-1} $ can be
holomorphically extended to $ \C^T $. \\ 
$ \bullet $ We chose $ U $ so that $ Z_U(Q,\mub,\cdot) $ is holomorphic on
$ \ds \left\{ \sb \in \C^T \ | \ \forall t \in \{1,\dots,T \} \ \si_t > -a
\right\} $. \\
So the formula $(*)$ gives the result. \\

\noindent {\bf Step 4:} $ R \in \er_{\ub} (Q) $ and $ S \in \er_{\vb} (R)
\Rightarrow S \in \er_{\ub+ \vb} (Q) $. \\

\noindent {\bf Proof of step 4}: \\
$ S $ is a linear combination of terms of the form $ \ds \p^{\bb} R
\prod_{t=1}^T \prod_{k \in F'_t} \p^{f'_t(k)} P_t $ where: \\
$ \bb \in \N^N $, $ F_1',\dots,F_T' $ are finite subset of $ \N $ such that
$ \forall t \ |F'_t|=v_t $, and $ f'_t \colon F'_t \to \N^N \setminus \{
\zerob \} $. \\
$ R \in  \er_{\ub} (Q) $ so $ \p^{\bb} R \in \er_{\ub} (Q) $ and then $
\p^{\bb} R $ is a linear combination of terms of the form: $ \ds  \p^{\gb}
Q \prod_{t=1}^T \prod_{k \in F_t} \p^{f_t(k)} P_t $ where: \\
$ \gb \in \N^N $, $ F_1,\dots,F_T $ are finite subset of $ \N $ such that $
\forall t \ |F_t|=u_t $, $ f_t \colon F_t \to \N^N \setminus \{ \zerob \}
$. \\

\noindent We can assume that $ \forall t,t' \in \{1,\dots,T \} \ \ F_t \cap
F'_{t'} = \emptyset $. \\

To conclude it is enough to show that: \\
 $ \ds U \overset{\text{def}}{=} \p^{\gb} Q \left( \prod_{t=1}^T \prod_{k
\in F_t} \p^{f_t(k)} P_t \right) \left( \prod_{t=1}^T \prod_{k \in F'_t}
\p^{f'_t(k)} P_t \right) $ is in $ \er_{\ub+ \vb} (Q) $. \\
For $ t \in \{1,\dots,T \} $ we define $ g_t \colon F_t \sqcup F'_t \to
\N^N \setminus \{ \zerob \} $ in the following way: \\
$ g_t(k)=f_t(k) $ if $ k \in F_t $, \ \  $ g_t(k)=f'_t(k) $ if $ k \in F'_t
$. \\
Then $ \ds U = \p^{\gb} Q \prod_{t=1}^T \prod_{k \in F_t\sqcup F'_t}
\p^{g_t(k)} P_t $ and $ \forall t \in \{1,\dots,T\} \ \ 
|F_t \sqcup F'_t| = u_t+v_t $. \\
So it is now clear that 
$ U \in \er_{\ub+ \vb} (Q) $.

\noindent {\bf Step 5}: \\
Let $ Q \in \R[X_1,\dots,X_N] $, $ a \in \R $ and $ b \in \N^* $. \\
Then $ Z(Q,\mub,\sb) $ is an entire combination until $ -a $ of functions
of the type $ Z(R,\mub,\sb+\ub) $ where $ \ub \in \N^N $ satisfies $ |\ub|
\geq b $ and $ R \in \er_{\ub}(Q) $. \\

\noindent {\bf Proof of step 5}: \\
the proof is by recurrence on $ b \in \N^* $. \\
For $ b = 1 $, it comes from step 3. \\
The combination of step 3 and step 4 allows us to deduce the result at rank
$ b+1 $ from the result at rank $ b $. \\

\noindent {\bf Last step: conclusion}: \\
let $ Q \in \R[X_1,\dots,X_N] $ and $ a \in \R_+ $. \\
We wish to show that $ Z(Q,\mub,\cdot) $ can be holomorphically extended
until $-a $. \\
Let $ b \in \N $. The value of $ b $ will be precised in the sequel. \\ 
By step 5 $ Z(Q,\mub,\sb) $ is an entire combination until $ -a $ of
functions of the type $ Z(R,\mub,\sb+ \ub) $ where $ \ub \in \N^N $
satisfies $ |\ub| \geq b $ and $ R \in \er_{\ub}(Q) $. \\
Let us consider $ \ub \in \N^N $ satisfying $ |\ub| \geq b $ and $ R \in
\er_{\ub}(Q) $. \\
$ R $ is a linear combination of polynomials of the form $ \ds S = \p^{\bb}
Q \prod_{t=1}^T \prod_{k \in F_t} \p^{f_t(k)} P_t $ where: \\ 
$ \bb \in \N^N $, $ F_1,\dots,F_T $ are finite subsets of $ \N $ satisfying
$ \forall t \ |F_t|=u_t $ and 
$ f_t \colon  F_t \to \N^N \setminus \{\zerob \} $. 

\begin{align*}
 \frac{\prod_{t=1}^T \prod_{k \in F_t} \p^{f_t(k)} P_t}{\prod_{t=1}^T
P_t^{u_t}}(\mb) & = \prod_{t=1}^T \prod_{k \in F_t} \frac{ \p^{f_t(k)}
P_t}{P_t}(\mb) \\ 
& \ll \prod_{t=1}^T  \prod_{k \in F_t} \mb^{- \ep \unb} \ \ (\mb \in
\N^{*N}) \\
& \ll \prod_{t=1}^T  \mb^{- \ep u_t \unb} \ \ (\mb \in \N^{*N}) \\
& \ll \mb^{- \ep |\ub| \unb} \ \ (\mb \in \N^{*N}) \\  
& \ll \mb^{- \ep b \unb}  \ \ (\mb \in \N^{*N})
\end{align*}
Let $ K $ be a compact of $ \C^T $ included in $ \ds \left\{ \sb \in \C^T \
| \ \forall t \in \{1,\dots,T \} \ \si_t > -a \right\} $. \\ 
As in step 2 : $ \ds \p^{\bb} Q(\mb) \prod_{t=1}^T P_t(\mb)^{-s_t} \ll
\mb^{q \unb + Tpa \unb} \ \ (\mb \in \N^{*N}\ \sb \in K) $. \\
$ \ds S(\mb) \prod_{t=1}^T P_t(\mb)^{-(s_t+u_t)} = \p^{\bb} Q(\mb)
\prod_{t=1}^T P_t(\mb)^{-s_t} \frac{\prod_{t=1}^T \prod_{k \in F_t}
\p^{f_t(k)} P_t}{\prod_{t=1}^T P_t^{u_t}}(\mb) $ \\
so $ \ds S(\mb) \prod_{t=1}^T P_t(\mb)^{-(s_t+u_t)} \ll \mb^{(q  + Tpa -
\ep b) \unb} \ \ (\mb \in \N^{*N}\ \sb \in K) $. \\
We choose $ b \in \N $ satisfying $ \ds b \geq \frac{q+Tpa+2}{\ep} $ so
that $ q  + Tpa - \ep b \leq -2 $. \\
We see that $ Z(R,\mub,\sb+\ub) $ is holomorphic on $ \ds \left\{ \sb \in
\C^T \ | \ \forall t \in \{1,\dots,T \} \ \si_t > -a \right\} $, so the
proof is done. \\
To obtain the formula of the theorem, it suffices to make $ \sb = -\kb $ in
the formula $(*)$ prooved in step 1 and to remark that $ Z_U(Q,\mub,-k) = 0 $
when we choose $ U $ such that $ U > \max \{k_1,\dots,k_T \} $. \\
This ends the proof of theorem 1. \\

\section{Proof of theorem 2}

Let $ \mub \in \left(\T \setminus \{ 1 \}\right)^N $ and $ \kb \in \N^{N}$ fixed.
Let $ d \in \N $ be fixed. \\
We denote $ \R_d[X_1,\dots,X_N] $ the set of the real polynomials of $ N $
variables with degree at most $ d $. \\
Let $ D = card \{ \ab \in \N^{N} \ | \ |\ab| \leq d \} $. \\
Let $ \phi \colon \R^D \to \ \R_d[X_1,\dots,X_N] $ defined by: 
$ \ds A =(a_{\ab})_{|\ab| \leq d} \mapsto \phi(A) = 
\sum_{|\ab| \leq d} a_{\ab}
\xb^{\ab} $. \\
It is an isomorphism of real vector spaces. \\
Thanks to theorem B of \cite{crise}, we know that there exists a polynomial
$ \ds G \in \R[X_1,\dots,X_{D(T+1)}] $ such that for all $ B,A_1,\dots,A_T
\in \R^D $ such that $ \phi(A_1),\dots, \phi(A_T)$ satisfy assumptions
of theorem II, we
have $Z(\phi(B); \phi(A_1),\dots, \phi(A_T); \mub;-\kb) =
G(B,A_1,\dots,A_T) $. \\

Theorem B of \cite{crise} implies that if we restrain to
polynomials of degree at most $ d $ with $ P_1,..,P_T $ satisfying
assumptions of theorem II,
the following relation, that we want to establish: \\ 
\begin{align*}
(1- \mub^{\lb}) Z(Q;P_1,\dots,P_T;\mub;-\kb) = & \mub^{\lb}  \sum_{ \zerob
< \ub \leq \kb }  \binom{\kb}{\ub}  Z \left( Q(\Xb + \lb) \prod_{t=1}^T
(\De P_t)^{u_t};P_1,\dots,P_T; \mub;-\kb + \ub \right) \\
                                   & + \mub^{\lb} Z(\De
Q;P_1,\dots,P_T;\mub;-\kb) + Z_{N-1}^{\lb}(-\kb)
\end{align*} 
is equivalent to $ G_1(B,A_1,\dots,A_T) =  G_2(B,A_1,\dots,A_T)  $ (with $
B = \phi^{-1}(Q) $ et $ \forall t \ A_t = \phi^{-1}(P_t) $) with $ G_1, G_2
\in \R[X_1,\dots,X_{D(T+1)}] $ depending only on $ \kb, d , N, T$ and $ \mub $
and repectively associated to the right side and left side. \\
 It is easy to see that if $ A \in \R_+^{*D} $ then $ \phi(A) $ is
hypoelliptic and non constant. \\
Therefore theorem I implies that for all $ B \in \R^D $, and for 
all $ A_1,\dots,A_T \in \R_+^{*D} $,
$ G_1(B,A_1,\dots,A_T) = G_2(B,A_1,\dots,A_T) $. Since $ G_1 $ and $ G_2 $
are polynomials, this implies that $ G_1 = G_2 $. \\ 
This end the poof of theorem 2. 

\section{Proof of corollaries:}
Corollary 1 is a direct consequence of theorem 1.\\
Point 1 of corollary 2 follows from assumption on ${\bf a}$ by easy 
computation.\\
{\bf Proof of the point 2 of corollary 2:}
By using theorem 1, to finish the proof of corollary 2, it's enough
to verifies that each polynomial $P_t$ is hypoelliptic.\\
Let $t\in \{1,\dots, T\}$ and $n\in \{1,\dots,N\}$ fixed. We have
uniformly in $[1,+\infty[^N$:
\begin{eqnarray*}
\frac{\partial P_t({\bf x})}{\partial x_n}
{P_t}^{-1}({\bf x})
&\ll & \frac{\sum_{k=1}^{r_t} \alpha_n^{t,k}|\langle \alpha^{t,k} ,
  {\bf x}\rangle | + c_{t,n} }
{\sum_{k=1}^{r_t} \left(\langle \alpha^{t,k} ,
  {\bf x}\rangle \right)^2 + \sum_{j=1}^N c_{t,j} x_j+d_t}\\
&\ll & \frac{1}
{\sqrt{\sum_{k=1}^{r_t} \left(\langle \alpha^{t,k} ,
  {\bf x}\rangle \right)^2 + \sum_{j=1}^N c_{t,j} x_j+d_t}}\\
&\ll& \frac{1}{(x_1+..+x_N)^{1/2}}.
\end{eqnarray*}
But deg$P_t \leq 2$. So the previous  implies that $P_t$ is an
hypoelliptic polynomial. This completes the proof of corollary 2.


\begin{thebibliography}{99}


\bibitem{akiyama; egami; tanigawa} {\bf Akiyama, Shigeki;  Egami, Shigeki;
Tanigawa, Yoshio} \\
Analytic continuation of multiple zeta-functions and their values at
non-positive integers. \\
Acta Arith. 98, No.2, 107-116 (2001). 


\bibitem{akiyama; ishikawa} {\bf Akiyama, Shigeki;  Ishikawa, Hideaki} \\
On analytic continuation of multiple $L$-functions and related zeta
functions. \\
Jia, C; K.Matsumoto (ed.) , Analytic number theory. Proceedings of the 1st
China-Japan seminar on number theory, Beijing, China, September 13-17, 1999
and the annual conference on analytic number theory, Kyoto, Japan, November
29-December 3, 1999. 
 Dordrecht: Kluwer Academic Publishers. Dev. Math. 6, 1-16 (2002)


\bibitem{akiyama; tanigawa} {\bf Akiyama, Shigeki; Tanigawa, Yoshio} \\
 Multiple zeta values at non-positive integers. 
The Ramanujan Journal, vol. 5, no.4 (2001) 327-351.


\bibitem{arakawa; kaneko} {\bf Arakawa, Tsuneo;  Kaneko, Masanobu} \\
Multiple zeta values, poly-Bernoulli numbers, and related zeta functions. \\
Nagoya Math. J. 153, 189-209 (1999).



\bibitem{cassou 79} {\bf Cassou-Noguès, Pierrette}. \\
Valeurs aux entiers négatifs des fonctions zêta et fonctions zêta
p-adiques. \\
Invent.Math. 51, 29-59 (1979). 


\bibitem{cassou 82}  {\bf Cassou-Noguès, Pierrette}. \\
Valeurs aux entiers négatifs de séries de Dirichlet associées à un
polyn\^ome.I. \\
J. Number Theory 14, 32-64 (1982).


\bibitem{chen; eie} {\bf  K.W. Chen; M. Eie}. \\
 A note on generalized Bernoulli numbers. 
 Pac. J. Math. 199, No.1, 41-59 (2001). 

\bibitem{crise} {\bf de Crisenoy, Marc}. \\
Values at T-tuples of negative integers of twisted multivariable zeta
series associated to polynomials of several variables. 
To appear in Compositio Mathematica. 


\bibitem{crise1} {\bf de Crisenoy, Marc}. \\
Valeurs aux $T-$uplets d'entiers négatifs de séries zêtas multivariables
associées à des polynômes de plusiurs variables. 
Thèse de l'université de Caen. (2003)

\bibitem{egami; matsumoto} {\bf Egami, Shigeki; Matsumoto, Kohji} \\
Asymptotic expansions of multiple zeta functions and power mean values of
Hurwitz zeta functions. 
J. Lond. Math. Soc., II. Ser. 66, No.1, 41-60 (2002). 

\bibitem{essouabri1} {\bf D. Essouabri} \\
Singularit\'es des s\'eries de Dirichlet associ\'ees \`a des polyn\^omes
de plusieurs variables et application \`a la th\'eorie analytique des
nombres.
Annales de l'institut Fourier, 47 (2), p.429-484 (1997).

\bibitem{essouabri} {\bf D. Essouabri} \\
Zeta function associated to Pascal's triangle mod $p$.\\
To appear in Japanese Journal of Mathematics.

\bibitem{hormander} {\bf L. H{\"o}rmander} \\
Analysis of linear partial differential operators II. Grundlehren,
vol. 257, Springer-Verlag, (1983).\\

\bibitem{lichtin} {\bf B. Lichtin} \\
On the moderate growth of generalized Dirichlet series for hypoelliptic polynomials.\\
Compositio Mathematica 80, No.3, 337-354 (1991).

\bibitem{shintani} {\bf T. Shintani} \\
 On evaluation of zeta functions of totally real algebraic number fields at
non-positive integers. \\
 J. Fac. Sci., Univ. Tokyo, Sect. I A 23, 393-417 (1976). 

\bibitem{zagier 94}{\bf Zagier, Don}. \\
Values of zeta functions and their applications. \\
Joseph, A. (ed.) et al., First European congress of mathematics (ECM),
Paris, France, July 6-10, 1992. Volume II: Invited lectures (Part 2).
Basel: Birkhäuser. Prog. Math. 120, 497-512 (1994). 


\bibitem{zhao} {\bf Zhao, Jianqiang} \\
Analytic continuation of multiple zeta functions. 
Proc. Am. Math. Soc. 128, No.5, 1275-1283 (2000). 

\end{thebibliography}
\end{document}